\documentclass[12pt]{article}
\usepackage[cp1251]{inputenc}
\usepackage{amssymb,latexsym,amsmath}

\usepackage{amsfonts}
\usepackage[english]{babel}
\usepackage{indentfirst}
\usepackage{graphicx,graphics,verbatim}
\usepackage{cite}
\usepackage{mathrsfs}

\newtheorem{theorem}{Theorem}

\begin{document}

\begin{center}		
	{\bf Inversion of the Weighted Spherical Mean} \\
\end{center}	

\begin{center}
	\textbf{Elina L. Shishkina}\\
	Voronezh State University, Voronezh, Russia\\
	shishkina@amm.vsu.ru\\
\end{center}


{\bf Abstract}. The paper contains the inversion formula for the weighted spherical mean. The interest to reconstruction a function by its integral by sphere grews tremendously in the last six decades, stimulated by the spectrum of new problems and methods of image reconstruction. We consider a generalization of the classical spherical mean and its inverse in the case when generalized translation acts to function instead of regular. As a particular case this problem includes  action of spherical means on radially symmetric functions.

\section{Introduction}
\label{sec:1}

Reconstruction of a function from a known subset of its spherical means is widely developed in pure and applied mathematics. 
Its connection with photoacoustic images is as follows.
Let the speed of sound propagation in the medium be a constant value. 
Then the pressure at a certain point in time is expressed in terms of the spherical mean pressure and its time derivative at some previous point in time \cite{Finch}.
Therefore, this imaging technique requires the inversion of spherical means.

The problem of reconstruction a function $f$ supported in a ball $B\in\mathbb{R}^n$, if the spherical means of $f$
are known over all geodesic spheres centered on the boundary $\partial B$ was solved using different approach in \cite{Finch,Kuchment,Rubin1,Rubin2,Rubin3,Agranovsky1,Agranovsky2,Kunyansky}.
It is remarkable that reconstruction formulas in \cite{Finch,Kuchment,Rubin1,Rubin2,Rubin3,Agranovsky1,Agranovsky2,Kunyansky} are different for even and odd dimension  of Euclidean space $n$.

Classical spherical mean has the form
\begin{equation}\label{ShSr}
	M(x,r,u)=\frac{1}{|S_{n}(1)|} \int\limits_{S_{n}(1)} u(x+\beta r)dS,\qquad x=(x_1,...,x_n),
\end{equation}
where $S_n(1)$ is unit sphere centered at the origin, $\beta$ is a coordinate of the sphere $S_n(1)$.

Operator (\ref{ShSr}) intertwines Laplace operator and  one-dimensional Bessel operator with index $n-1$. 
Great interest among various researchers is a generalization of the  spherical mean (\ref{ShSr}).
So, in paper \cite{Weinstein2} was considered a spherical mean in space with negative curvature, in \cite{Dunkl} and \cite{Elouadih} was studied a generalization of the spherical mean generated by the  Dunkl  transmutation operator.
In this paper we consider the spherical weighted mean (see  (\ref{05})), which is the transmutation operator intertwining the multidimensional operator 
\begin{equation}\label{LB}
	(\Delta_\gamma)_x{=}\sum\limits_{i=1}^n (B_{\gamma_i})_{x_i},\qquad (B_{\gamma_i})_{x_i}{=}\frac{\partial^2}{\partial x_i^2}+\frac{\gamma_i}{x_i}\frac{\partial}{\partial x_i},\qquad \gamma_i\geq0,\qquad i=1,...,n
\end{equation}
and one-dimensional Bessel operator with index $n+|\gamma|-1$ of the form
\begin{equation}\label{Bess}
	(B_{n+|\gamma|-1})_t{=}\frac{d^2}{d t^2}+\frac{n+|\gamma|-1}{t}\frac{d}{dt},\qquad t>0,\qquad |\gamma|{=}\gamma_1{+}{\ldots}{+}\gamma_{n}.
\end{equation}
Such spherical mean is closely related to B-ultra-hyperbolic equation of the form (see \cite{LPSh1,LPSh2})
\begin{equation}\label{UG22}
	\sum\limits_{j=1}^{n}(B_{\gamma_i})_{x_j}u=\sum\limits_{j=1}^{n}(B_{\gamma_i})_{y_j}u,\qquad u=u(x_1,...,x_n,y_1,...,y_n).
\end{equation}

Not so long ago, various new methods for solving inverse problems with the Bessel operator appeared.
The inverse problem involving recovery of initial temperature from the information of final temperature
profile in the case of heat equation with Bessel operator was studied in \cite{Masood}.
The  solution of the inverse spectral problem for Bessel-type differential operators
on noncompact star-type graphs under generalized Neumann-type boundary conditions was considered in \cite{Yurko}.
V. V. Kravchenko with coauthors introduced a new method to solution of the inverse Sturm-Liouville 
problem  (see \cite{Krav1,Krav2,KKK QuantumFest,Kr2019JIIP}) 
which can be adapted to inverse problem with Bessel operator.
Their idea is based on the observation that the
potential can be recovered from the very first coefficient of the
Fourier--Legendre series, and to find this coefficient a system of linear
algebraic equations can be obtained directly from the Gel'fand-Levitan
equation.

The paper contains an
inversion formula for the weighed spherical mean (\ref{05}) using the properties of the mixed Riesz hyperbolic B-potential (\ref{RPot1}).

\section{Basic Definitions}
\label{sec:2}

In this section we  give a summary of the basic notations, terminology and results which will be used
in this article.

Suppose that $\mathbb{R}^{n+1}$ is the $n+1$--dimensional Euclidean space,
$$
\mathbb{R}^{n+1}_+{=}\{(t,x){=}(t,x_1,\ldots,x_n)\in\mathbb{R}^{n+1},\,\,\, x_1{>}0,\ldots, x_n{>}0\},
$$
$\gamma{=}(\gamma_1,{...},\gamma_{n})$ is a multiindex  consisting of  fixed real numbers $\gamma_i\geq 0$, $i{=}1,{...},n,$   and $|\gamma|{=}\gamma_1{+}{\ldots}{+}\gamma_{n}$.
Let $\Omega$ be a finite or infinite open set in $\mathbb{R}^{n+1}$ symmetric with respect  to each hyperplane $x_i{=}0$, $i=1,...,n$, $\Omega_+=\Omega\cap{\mathbb{R}}^{n+1}_+$ and $\widetilde{\Omega}_+=\Omega\cap\overline{\mathbb{R}}\,\!^{n+1}_+$ where $$
\overline{\mathbb{R}}\,^{n+1}_+{=}\{(t,x){=}(t,x_1,\ldots,x_n){\in}\mathbb{R}^{n+1}, x_1{\geq}0,\ldots, x_n{\geq}0\}.
$$
We deal with the class $C^m(\Omega_+)$ consisting of $m$ times differentiable on $\Omega_+$ functions
and denote by $C^m(\widetilde{\Omega}_+)$ the subset of functions from $C^m(\Omega_+)$ such that all existing derivatives  of these functions with respect to $x_i$ for any $i=1,...,n$  are continuous up to $x_i{=}0$ and all existing derivative  with respect to $t$   are continuous for $t\in\mathbb{R}$.
Class $C^m_{ev}(\widetilde{\Omega}_+)$ consists of all functions from $C^m(\widetilde{\Omega}_+)$ such that $\frac{\partial^{2k+1}f}{\partial x_i^{2k+1}}\biggr|_{x=0}=0$ for all nonnegative integer $k\leq \frac{m-1}{2}$ and for $i=1,...,n$ (see \cite{Zhit} and \cite{kipr}, p. 21).
In the following we will denote $C^m_{ev}(\overline{\mathbb{R}}\,\!^{n+1}_+)$ by $C^m_{ev}$.
We set  
$$
C^\infty_{ev}(\widetilde{\Omega}_+)=\bigcap C^m_{ev}(\widetilde{\Omega}_+)
$$
with intersection taken for all finite $m$. Let $C^\infty_{ev}(\overline{\mathbb{R}}\,^{n+1}_+)=C^\infty_{ev}$.
Assuming that ${\stackrel{\circ}C}\,\!^\infty_{ev}(\widetilde{\Omega}_+)$ is the space of all functions  $f{\in}C^\infty_{ev}(\widetilde{\Omega}_+)$ with a compact support. We will use the notation ${\stackrel{\circ}C}\,\!^\infty_{ev}(\widetilde{\Omega}_+){=}\mathcal{D}_+(\widetilde{\Omega}_+)$.

Let $\mathcal{L}_p^{\gamma}(\Omega_+)$, $1{\leq}p{<}\infty$ be the space of  all measurable in $\Omega_+$ functions such that
$$
\int\limits_{\Omega_+}|f(t,x)|^p x^\gamma dtdx<\infty,
$$
where and further
$$
x^\gamma=\prod\limits_{i=1}^n x_i^{\gamma_i}.
$$
For a real number $p\geq 1$, the $\mathcal{L}_p^\gamma(\Omega_+)$--norm of $f$ is defined by
$$
||f||_{\mathcal{L}_p^\gamma(\Omega_+)}=\left(\,\,\int\limits_{\Omega_+}|f(t,x)|^p x^\gamma dtdx\right)^{1/p}.
$$
Let $\mathcal{L}_p^\gamma=\mathcal{L}_p^\gamma(\mathbb{R}^+_{n+1})$.

We will use the generalized convolution product defined  by the formula
\begin{equation}\label{Conv}
	(f*g)_\gamma(x,t)=\int\limits_{\mathbb{R}^{n+1}_+}f(\tau,y)(\,^\gamma \mathbb{T}^y_xg)(t-\tau,x)y^\gamma d\tau dy,
\end{equation}
where $^\gamma \mathbb{T}^y_x$ is multidimensional generalized translation
\begin{equation}\label{Transl}
	(^\gamma \mathbb{T}^y_xf)(t,x){=}(^{\gamma_1}T_{x_1}^{y_1}{...}^{\gamma_n}T_{x_n}^{y_n}f)(t,x).
\end{equation}
Each of  one--dimensional generalized translations $^{\gamma_i}T_{x_i}^{y_i}$ is defined for $i{=}1,{...},n$
by the next formula (see \cite{Lev1}, p. 122, formula (5.19))
$$
(^{\gamma_i}T_{x_i}^{y_i}f)(t,x){=}
\frac{\Gamma\left(\frac{\gamma_i+1}{2}\right)}{\Gamma\left(\frac{\gamma_i}{2}\right)\Gamma\left(\frac{1}{2}\right)}\int\limits_0^\pi\sin^{\gamma_i-1}\varphi_i \times$$
$$\times
f(t,x_1,...,x_{i-1},\sqrt{x _i^2+y_i^2-2x_iy_i\cos\varphi_i},x_{i+1},...,x_{n})\,\,
d\varphi_i,
$$
$\gamma_i>0,$ $i=1,...,n$ and
for $\gamma_i=0$ generalized translation $\,^{\gamma_i} T^{y_i}_{x_i}$ is
$$
\,^0 T^{y_i}_{x_i}= \frac{f(x+y)-f(x-y)}{2}.
$$

As the space of basic functions we will use  the subspace  of rapidly decreasing functions:
$$
S_{ev}({\mathbb{R}}^{n+1}_+)=\left\{f\in C^\infty_{ev}:\sup _{{(t,x)\in {{R}}^{n+1}_+}}\left|t^{\alpha_0}x^{\alpha }D^{\beta }f(t,x)\right|<\infty \right\},
$$
where $\alpha=(\alpha_1,...,\alpha_n)$, $\beta=(\beta_0,\beta_1,...,\beta_n)$, $\alpha_0,\alpha_1,...,\alpha_n,\beta_0,\beta_1,...,\beta_n$ are arbitrary integer nonnegative numbers, $x^\alpha{=}x_1^{\alpha_1} x_2^{\alpha_2} \ldots x_n^{\alpha_n}$,
${D}^\beta{=}{D}^{\beta_0}_{t}{D}^{\beta_1}_{x_1}{...}{D}^{\beta_n}_{x_n}$, ${D}_{t}=\frac{\partial}{\partial t}$, ${D}_{x_j}=\frac{\partial}{\partial x_j}$, $j=1,...,n$.

The multidimensional Fourier--Bessel transform  of  a function $f{\in} \mathcal{L}^\gamma_1(\mathbb{R}^{n+1}_+)$ is
$$
\mathcal{F}_\gamma[f](\tau,\xi)=\widehat{f}(\tau,\xi)=\int\limits_{\mathbb{R}^{n+1}_+}f(t,x)\,e^{-it\tau}\,\mathbf{j}_\gamma(x;\xi)x^\gamma dtdx,
$$
where
$$\mathbf{j}_\gamma(x;\xi)=\prod\limits_{i=1}^n j_{\frac{\gamma_i-1}{2}}(x_i\xi_i),\qquad \gamma_1\geq0,...,\gamma_n\geq0.$$

Now we introduce a  weighted spherical mean.
When constructing a  weighted spherical mean, instead of the usual shift, a multidimensional generalized translation (\ref{Transl}) is used.

Weighted spherical mean (see \cite{LPSh1,LPSh2,SitnikShishkina Elsevier}) of function $f(x)$, $x\in \overline{\mathbb{R}}\,^n_+$ for $n\geq 2$ is 
\begin{equation}\label{05}
	(M^\gamma_tf)(x)=(M^\gamma_t)_x[f(x)]=\frac{1}{|S_1^+(n)|_\gamma}\int\limits_{S^+_1(n)}\,^\gamma
	\mathbb{T}_x^{t\theta}f(x)
	\theta^\gamma dS,
\end{equation}
where $\theta^\gamma{=}\prod\limits_{i=1}^{n}\theta_i^{\gamma_i},$ $S^+_1(n){=}\{\theta{:}|\theta|{=}1,\theta{\in}\mathbb{R}^n_+\}$ is a part of a sphere in
$\mathbb{R}^n_+$, and $|S^+_1(n)|_\gamma$ is given by 
\begin{equation}
	\label{150}
|S_1^+(n)|_\gamma
=\int\limits_{S^+_1(n)} x^\gamma dS=
\frac{\prod\limits_{i=1}^n{\Gamma\left(\frac{\gamma_i{+}1}{2}\right)}}{2^{n-1}\Gamma\left(\frac{n{+}|\gamma|}{2}\right)}.
\end{equation}
For $n=1$ let $M^\gamma_t[f(x)]=\,^\gamma {T}_x^{t}f(x)$.

Let $\nu> 0$.	One-dimensional Poisson operator is defined for integrable  function $f$ by the equality
\begin{equation}
	\label{154}
	\mathcal{P}_x^{\nu}f(x)=\frac{2C(\nu)}{x^{\nu-1}}
	\int\limits_0^x \left( x^2-t^2\right)^{\frac{\nu}{2}-1}f(t)\,dt,\qquad   C(\nu)=	\frac{\Gamma\left(\frac{\nu+1}{2}\right)}{\sqrt{\pi}\,\Gamma\left(\frac{\nu}{2}\right)},
\end{equation}
or
\begin{equation}
	\label{155}
	\mathcal{P}_x^{\nu}f(x)=
	C(\nu)
	\int\limits_0^\pi f(x\cos\varphi)\sin^{\nu-1}\varphi\:d\varphi,.
\end{equation}
The constant $C(\nu)$ is chosen so that $\mathcal{P}_x^\nu[1]=1$.
For $\nu=0$ one-dimensional Poisson operator turns into an identical operator: $\mathcal{P}_x^{0}=I$.

\section{Mixed Hyperbolic Riesz B--potential and Its Inversion}
\label{sec:3}

The given supplementary information in this section will be used in the 4th section where the main results are presented.

Mixed hyperbolic Riesz B--potential was studied in \cite{ShishkinaKir,ShishkinaRUDN,SitnikShishkina Elsevier}. 
Here we provide  a definition of the mixed hyperbolic Riesz B--potential and give some separation results that we will use in the next section. 

Let $|x|=\sqrt{x_1^2+...+x_n^2}$. First for $(t,x)\in\mathbb{R}^{n+1}_+$, $\lambda\in{C}$ we define  function $s^\lambda$  by the formula
\begin{equation}\label{FunS}
	s^\lambda(t,x)=\left\{
	\begin{array}{ll}
		$$\frac{(t^2-|x|^2)^\lambda}{N(\alpha,\gamma,n)}$$, & \hbox{{when}\, $t^2\geq |x|^2$\, {and}\, $t\geq 0$;} \\
		$$0,$$ & \hbox{{when}\, $t^2<|x|^2$\, {or}\, $t<0$,}
	\end{array}
	\right.
\end{equation}
where
\begin{equation}\label{Const}
	N(\alpha,\gamma,n)=\frac{2^{\alpha-n-1}}{\sqrt{\pi}}\prod\limits^{n}_{i=1}\Gamma\left(\frac{\gamma_i+1}{2}\right)
	\Gamma\left(\frac{\alpha-n-|\gamma|+1}{2}\right)\Gamma\left(\frac{\alpha}{2}\right).
\end{equation}
Regular weighted distribution corresponding
to (\ref{FunS}) we will denote by $s_+^\lambda$.

We  introduce the \emph{mixed  hyperbolic Riesz B--potential $I_{s_,\gamma}^\alpha$ of order} $\alpha$ as a
generalized convolution product (\ref{Conv})  with  a weighted distribution $s_+^{\frac{\alpha-n-|\gamma|-1}{2}}$ and $f\in S_{ev}$ (see \cite{ShishkinaRUDN}):
\begin{equation}\label{RPot1}
	(I_{s,\gamma}^\alpha
	f)(t,x)=\left( s_+^{\frac{\alpha-n-|\gamma|-1}{2}}*f\right)_\gamma(t,x).
\end{equation}
The precise definition of the constant $N(\alpha,\gamma,n)$ allows  to obtain the semigroup property or index low of the potential (\ref{RPot1}).

We can rewrite formula (\ref{RPot1}) as
\begin{equation}\label{RPot2}
	(I_{s,\gamma}^\alpha
	f)(t,x)=\int\limits_{\mathbb{R}^{n+1}_+}s^{\frac{\alpha-n-|\gamma|-1}{2}}_+(\tau,y)(\,^{\gamma}\mathbb{T}^{y}_{x})f(t-\tau,x)y^\gamma
	d\tau dy.
\end{equation}
Integral (\ref{RPot2}) converges absolutely for $n{+}|\gamma|{-}1{<}\alpha$ for integrable with weight $y^\gamma$ on the part of the cone $\{|y|<\tau\}^+{=}\{y{\in}\mathbb{R}^n_+{:}|y|{<}\tau\}$, $0<\tau<t$
function $f(\tau,y)$.

For $0\leq \alpha \leq n+|\gamma|-1$
$$
(I_{s,\gamma}^\alpha
f)(t,x)= \left(\frac{\partial^2}{\partial t^2}-\Delta_\gamma \right)^q	(I_{s,\gamma}^{\alpha+2q}
f)(t,x)
$$
where $q=\left[\frac{n+|\gamma|-\alpha+1}{2} \right] $.

In the case when $f(t,x)=h(t)F(x)$ we get
$$
	(I_{s,\gamma}^{\alpha}
	hF)(t,x)=
	$$
	\begin{equation}\label{RPot3}
	=\frac{1}{N(\alpha,\gamma,n)}\int\limits_0^\infty h(t-\tau)d\tau \int\limits_{|y|<\tau}(\tau^2-|y|^2)^{\frac{\alpha-n-|\gamma|-1}{2}}(\,^{\gamma}\mathbb{T}^{y}_{x})F(x)y^\gamma dy.
\end{equation}

\begin{theorem}\cite{ShishkinaKir}
	Let $n+|\gamma|-1{<}\alpha{<}n+|\gamma|+1$,
	$1\leq p{<}\frac{n+|\gamma|+1}{\alpha}$. For the next estimate
	\begin{equation}\label{Eq4}
		||I_{s,\gamma}^\alpha
		f||_{q,\gamma} \leq M||f||_{p,\gamma},\qquad f\in
		S_{ev}
	\end{equation}
	to be valid it is necessary and sufficient that
	$q=\frac{(n+|\gamma|+1)p}{n+|\gamma|+1-\alpha p}$. Constant $M$ does not depend on $f$.
\end{theorem}

\textbf{Remark.} By virtue of (\ref{Eq4}) there is unique extension of $I_{s,\gamma}^\alpha$ to all $\mathcal{L}_p^\gamma$, $1<p<\frac{n+|\gamma|+1}{\alpha}$ preserving boundedness when $n+|\gamma|-1<\alpha<n+|\gamma|$.
It follows that this extension is introduced by the integral (\ref{RPot2})
from its absolute convergence.

\begin{theorem}\cite{ShishkinaKir}
		For $f\in S_{ev}$  the Fourier--Hankel transform of mixed hyperbolic Riesz potential $I_{s,\gamma}^\alpha
	f$ is 
	\begin{equation}\label{HankelPot}
		\mathcal{F}_\gamma[I_{s,\gamma}^\alpha
		f](\tau,\xi)=\mathscr{Q}\left|\tau^2-|\xi|^2\right|^{-\frac{\alpha}{2}}\cdot\mathcal{F}_\gamma[f(t,x)](\tau,\xi),
	\end{equation}
	where
	$$
	\mathscr{Q}=\left\{
	\begin{array}{ll}
	$$1$$, & \hbox{$|\xi|^2\geq \tau^2$;} \\
	$$e^{-\frac{\alpha\pi}{2}i}$$, & \hbox{$|\xi|^2<\tau^2$, $\tau\geq0$;} \\
	$$e^{\frac{\alpha\pi}{2}i}$$, & \hbox{$|\xi|^2<\tau^2$, $\tau<0$.}
	\end{array}
	\right.
	$$
\end{theorem}

For the inversion of the potential (\ref{RPot1}) 
approach  based on the idea of approximative inverse operators (see \cite{Nogin3}) was used.  This method gives an inverse operator as a limit of regularized operators.
Namely,  taking into account the formula (\ref{HankelPot}) we will construct inverse operator for the  potential (\ref{RPot1})   in the form
$$
(I_{s,\gamma}^\alpha)^{-1}
f=\lim\limits_{\varepsilon\rightarrow 0}\left(\mathcal{F}_\gamma^{-1}(\mathscr{Q}|\tau^2-|\xi|^2|^{\frac{\alpha}{2}}e^{-\varepsilon|\tau|-\varepsilon|\xi|})\ast f\right)_\gamma,
$$
where the limit is understood in the norm $\mathcal{L}_p^\gamma$ or almost everywhere.

Let
$$
g_{\alpha,\gamma,\varepsilon}(t,x)=\mathcal{F}_\gamma^{-1}(\mathscr{Q}|\tau^2-|\xi|^2|^{\frac{\alpha}{2}}e^{-\varepsilon|\tau|-\varepsilon|\xi|})(t,x),
$$
then
\begin{equation}\label{InvOp}
	(I_{s,\gamma}^\alpha)^{-1}
	f=\lim\limits_{\varepsilon\rightarrow 0}(g_{\alpha,\gamma,\varepsilon}*f)_\gamma.
\end{equation}

\begin{theorem}\cite{ShishkinaKir}
	The function $g_{\alpha,\gamma,\varepsilon}(t,x)$ belongs to the space $\mathcal{L}_p^\gamma,$ $1<p<\infty$ with additional restriction $\frac{2(n+|\gamma|)-1}{2(n+|\gamma|)-2}<p$ for $n+|\gamma|-1<\alpha<n+|\gamma|$ when $n+|\gamma|+1$ is odd.
\end{theorem}

\begin{theorem}\label{Inv}\cite{ShishkinaKir}
	Let $n+|\gamma|-1<\alpha<n+1+|\gamma|$, $1<p<\frac{n+1+|\gamma|}{\alpha}$ with the additional restriction $p<\frac{2(n+1+|\gamma|)(n+|\gamma|)}{n+1+|\gamma|+2\alpha(n+|\gamma|)}$ when $n+|\gamma|-1<\alpha<n+|\gamma|$ and $n$ is odd. Then
	$$
	((I_{s,\gamma}^\alpha)^{-1}I_{s,\gamma}^\alpha f)(t,x)= f(t,x),\qquad f(t,x)\in \mathcal{L}_p^\gamma,
	$$
	where $(I_{s,\gamma}^\alpha)^{-1}f=\lim\limits_{\varepsilon\rightarrow0}(I_{s,\gamma}^\alpha)^{-1}_\varepsilon f$.
\end{theorem}

It is easy to see that if $\alpha=2m$, $m\in\mathbb{N}$ the inverse to $I_{s,\gamma}^{2m}$
 is $\left(\frac{\partial^2}{\partial t^2}-\Delta_\gamma \right)^m$:
 $$
 \left(\frac{\partial^2}{\partial t^2}-\Delta_\gamma \right)^m(I_{s,\gamma}^{2m} f)(t,x)=f(t,x),\qquad f(t,x)\in \mathcal{L}_p^\gamma.
 $$

  The inverse  Fourier--Hankel transform of $\mathscr{Q}|\tau^2-|\xi|^2|^{\frac{\alpha}{2}}e^{-\varepsilon|\tau|-\varepsilon|\xi|}$ can be represented in the form
	$$
	g_{\alpha,\gamma,\varepsilon}(t,x)=\mathscr{C}(n,\gamma,\alpha)\int\limits_0^\infty |1-r^2|^{\frac{\alpha}{2}}\,r^{n+|\gamma|-1}\times
	$$	
	\begin{equation}\label{Yadro}
		\times \left[ e^{-\frac{\alpha\pi i}{2}\theta(1-r)}\mathscr{F}_{n,\gamma,\alpha,\varepsilon}^+(r,x,t)+ e^{\frac{\alpha\pi i}{2}\theta(1-r)}\mathscr{F}_{n,\gamma,\alpha,\varepsilon}^-(r,x,t) \right]dr,
	\end{equation}
	where
	$$
	\mathscr{C}(n,\gamma,\alpha)=\frac{\prod\limits^n_{i=1}\Gamma\left(\frac{\gamma_i+1}{2}\right)}{2^{n-1}}
	\frac{\Gamma(n+|\gamma|+1+\alpha)}{\Gamma\left(\frac{n+|\gamma|}{2}\right)}
	$$
	$$
	\mathscr{F}_{n,\gamma,\alpha,\varepsilon}^+(r,x,t)= \frac{\,_2F_1\left(\frac{n+|\gamma|+1+\alpha}{2},\frac{n+|\gamma|+2+\alpha}{2};\frac{n+|\gamma|}{2};-\frac{|x|^2r^2}{(\varepsilon +\varepsilon r+it)^2}\right)}{{(\varepsilon +\varepsilon r+it)^{n+|\gamma|+1+\alpha}}},
	$$
	$$
	\mathscr{F}_{n,\gamma,\alpha,\varepsilon}^-(r,x,t)= \frac{\,_2F_1\left(\frac{n+|\gamma|+1+\alpha}{2},\frac{n+|\gamma|+2+\alpha}{2};\frac{n+|\gamma|}{2};-\frac{|x|^2r^2}{(\varepsilon +\varepsilon r-it)^2}\right)}{{(\varepsilon +\varepsilon r-it)^{n+|\gamma|+1+\alpha}}}.
	$$

\section{Inverse Problem}

In this section we
consider the recovery of a  function $f$ from a knowledge of its weighted spherical mean $M^\gamma_\rho f$.

Let $f=f(x)\in C^2(\mathbb{R}^n_+)$, such that $\frac{\partial f}{\partial x_i}\biggr|_{x_i=0}=0$, $i=1,...,n$.
The weighted spherical mean $M^\gamma_t f$ is the transmutation operator intertwining $(\Delta_\gamma)_x$ and $(B_{n+|\gamma|-1})_t$ (see \cite{ShishkinaSitnik}):
\begin{equation}\label{OPPr}
	(B_{n+|\gamma|-1})_t(M^\gamma_tf)(x)=(M^\gamma_t(\Delta_\gamma)_xf)(x).
\end{equation}

Let consider the integral operator
\begin{equation}\label{IntOp}
	(\mathscr{M}^{\gamma,k}_tf)(x)=\frac{1}{|S_1^+(n)|_\gamma}\int\limits_{\mathbb{R}^n_+}  (\,^\gamma\mathbf{T}_x^y f(x))(t^2-|y|^2)_+^{\frac{k-n-|\gamma|-1}{2}}y^\gamma dy.
\end{equation}
Operator $t^{1-k}\mathscr{M}^{\gamma,k}_t$ intertwines   $(\Delta_\gamma)_x$ and $(B_{k})_t$ when $k>n+|\gamma|-1$:
\begin{equation}\label{OPPr2}
	(B_{k})_t(t^{1-k}\mathscr{M}^{\gamma,k}_tf)(x)=(t^{1-k}\mathscr{M}^{\gamma,k}_t(\Delta_\gamma)_xf)(x).
\end{equation}

We'll tend to use spherical coordinates in (\ref{IntOp}) when $k>n+|\gamma|-1$, then, using (\ref{154}) we can write
$$
(\mathscr{M}^{\gamma,k}_tf)(x)
=
$$
$$
=\frac{1}{|S_1^+(n)|_\gamma}\int\limits_{\{|y|<t\}^+}  (\,^\gamma\mathbf{T}_x^y f(x))(t^2-|y|^2)^{\frac{k-n-|\gamma|-1}{2}}y^\gamma dy=\{y=\rho\theta\}=
$$
$$
=\frac{1}{|S_1^+(n)|_\gamma}\int\limits_0^t (t^2-\rho^2)^{\frac{k-n-|\gamma|-1}{2}}\rho^{n+|\gamma|-1}d\rho \int\limits_{S^+_1(n)} (\,^\gamma\mathbf{T}_x^{\rho\theta} f(x))\theta^\gamma dS=
$$
$$
=\int\limits_0^t (t^2-\rho^2)^{\frac{k-n-|\gamma|-1}{2}}\rho^{n+|\gamma|-1}(M^\gamma_\rho f)(x)d\rho=
$$
$$
=\frac{t^{k-n-|\gamma|}}{2C(k-n-|\gamma|+1)} \left( \mathcal{P}_t^{k-n-|\gamma|+1}t^{n+|\gamma|-1}(M^\gamma_\rho f)(x)\right) (t).
$$

Now let find the inverse operator for $\mathscr{M}^{\gamma,k}$.
Let multiply (\ref{IntOp}) by $ h(t-\tau)$ and integrate by $\tau$ from $0$ to $\infty$. 
The function $h(t)$ should be chosen such that the function $h(t-\tau)(\mathscr{M}^{\gamma,k}_\tau f)(x)$ is an integrable by $\tau$ by the interval from $0$ to $\infty$. 

We obtain
$$
\int\limits_0^\infty h(t-\tau)(\mathscr{M}^{\gamma,k}_\tau f)(x)d\tau=
$$
$$
=\frac{1}{|S_1^+(n)|_\gamma}\int\limits_0^\infty h(t-\tau)d\tau
\int\limits_{\{|y|<\tau\}^+}  (\tau^2-|y|^2)^{\frac{k-n-|\gamma|-1}{2}}(\,^\gamma\mathbf{T}_x^y f)(x)y^\gamma dy.
$$
Taking into account (\ref{RPot3}) we get
$$
\frac{|S_1^+(n)|_\gamma}{N(k,\gamma,n)}\int\limits_0^\infty h(t-\tau)(\mathscr{M}^{\gamma,k}_\tau f)(x)d\tau=(I_{s,\gamma}^{k}
hf)(t,x),
$$
where $N(k,\gamma,n)$ defined by
(\ref{Const}) and $I_{s,\gamma}^k$ is  the mixed  hyperbolic Riesz B--potential  (\ref{RPot2}) of order $k>0$ acting to function $h(t)f(x)$.
Therefore, using theorem \ref{Inv}, we obtain
\begin{equation}\label{inv}
	h(t)f(x)=\frac{|S_1^+(n)|_\gamma}{N(k,\gamma,n)}\left( (I_{s,\gamma}^k)^{-1}_{\phi,y} \int\limits_0^\infty h(\phi-\tau)(\mathscr{M}^{\gamma,k}_\tau f)(y)d\tau\right)(t,x), 
\end{equation}
where $n+|\gamma|-1<k<n+1+|\gamma|$ and
$
(I_{s,\gamma}^k)^{-1}
$ is given by (\ref{InvOp}).
So, in the inverse formula (\ref{inv}) we have an arbitrary parameter $k\in(n+|\gamma|-1,n+1+|\gamma|)$ and   an arbitrary non-zero function $h$ (such that the function $h(t-\tau)(\mathscr{M}^{\gamma,k}_\tau f)(x)$ is an integrable by $\tau$) depending on one variable.

In order to find the inverse to weighed spherical mean the formula (\ref{inv}) can be simplified. We can take 
$k=2m>n+|\gamma|-1$, $m\in\mathbb{N}$. In this case $$(I_{s,\gamma}^{2m})^{-1}=\left(\frac{\partial^2}{\partial t^2}-\Delta_\gamma \right)^m $$
and
\begin{equation}\label{inv1}
	h(t)f(x)=\frac{|S_1^+(n)|_\gamma}{N(2m,\gamma,n)}\left(\frac{\partial^2}{\partial t^2}-\Delta_\gamma \right)^m  \int\limits_0^\infty h(t-\tau)(\mathscr{M}^{\gamma,2m}_\tau f)(x)d\tau. 
\end{equation}

So we obtain the main statement.

\begin{theorem}
Let $f=f(x)\in C^2(\mathbb{R}^n_+)$, such that $\frac{\partial f}{\partial x_i}\biggr|_{x_i=0}=0$, $i=1,...,n$ and 
$$
(\mathscr{M}^{\gamma,k}_tf)(x)
=\frac{t^{k-n-|\gamma|}}{2C(k-n-|\gamma|+1)} \left( \mathcal{P}_t^{k-n-|\gamma|+1}t^{n+|\gamma|-1}(M^\gamma_\rho f)(x)\right) (t),
$$
where $M^\gamma_\rho f$ is the weighted spherical mean (\ref{05}) of the function $f$, 
is $\mathcal{P}_t^{\nu}$ the one-dimensional Poisson operator (\ref{154}), $C(\nu)$ is the constant defined in (\ref{154}). Then the function $f$
can be reconstructed by its  weighted spherical mean by the formula
$$
	h(t)f(x)=\frac{|S_1^+(n)|_\gamma}{N(2m,\gamma,n)}\left(\frac{\partial^2}{\partial t^2}-\Delta_\gamma \right)^m  \int\limits_0^\infty h(t-\tau)(\mathscr{M}^{\gamma,2m}_\tau f)(x)d\tau, 
$$
where  function $h(t)$ is arbitrary such that the function $h(t-\tau)(\mathscr{M}^{\gamma,k}_\tau f)(x)$ is an integrable by $\tau$ by the interval from $0$ to $\infty$, $|S_1^+(n)|_\gamma$ is given by (\ref{150}), $N(2m,\gamma,n)$ is given by (\ref{Const}).
\end{theorem}

{\bf Example.} Let $h(t)=e^t$, $(M^\gamma_\rho f)(x)={\bf j}_\gamma(x,\xi)\,\,j_{\frac{n+|\gamma|}{2}-1}(\rho|\xi|),$
where $\xi=(\xi_1,...,\xi_n)$ some vector.
$$
(\mathscr{M}^{\gamma,2m}_\tau f)(x)
=\int\limits_0^\tau (\tau^2-\rho^2)^{\frac{2m-n-|\gamma|-1}{2}}\rho^{n+|\gamma|-1}(M^\gamma_\rho f)(x)d\rho=
$$
$$
={\bf j}_\gamma(x,\xi)\,\int\limits_0^\tau (\tau^2-\rho^2)^{\frac{2m-n-|\gamma|-1}{2}}\rho^{n+|\gamma|-1}\,j_{\frac{n+|\gamma|}{2}-1}(\rho|\xi|)d\rho=
$$
$$
=2^{\frac{n+|\gamma|}{2}-1}\,\Gamma\left(\frac{n+|\gamma|}{2} \right)|\xi|^{1-\frac{n+|\gamma|}{2}}{\bf j}_\gamma(x,\xi)\times
$$
$$
\times\int\limits_0^\tau \left( \tau^2-\rho^2\right)^{\frac{2m-n-|\gamma|-1}{2}}\rho^{\frac{n+|\gamma|}{2}}J_{\frac{n+|\gamma|}{2}-1}(\rho|\xi|)\,d\rho.
$$
Using formula 2.12.4.6 from \cite{IR2} we obtain
$$
(\mathscr{M}^{\gamma,2m}_\tau f)(x)
=
$$
$$
=2^{\frac{n+|\gamma|}{2}-1}\,\Gamma\left(\frac{n+|\gamma|}{2} \right)|\xi|^{1-\frac{n+|\gamma|}{2}}{\bf j}_\gamma(x,\xi)\frac{2^{\frac{2m-n-|\gamma|-1}{2}}\tau^{m-\frac{1}{2}}}{|\xi|^{\frac{2m-n-|\gamma|+1}{2}}}\times
$$
$$
\times\Gamma\left(\frac{2m-n-|\gamma|+1}{2} \right) J_{m-\frac{1}{2}}(|\xi|\tau)=
$$
$$
=\frac{2^{m-\frac{3}{2}}\tau^{m-\frac{1}{2}}}{|\xi|^{m-\frac{1}{2}}}\,\Gamma\left(\frac{n+|\gamma|}{2} \right)\Gamma\left(\frac{2m-n-|\gamma|+1}{2} \right) {\bf j}_\gamma(x,\xi)J_{m-\frac{1}{2}}(|\xi|\tau).
$$
Taking into account that $h(t)=e^t$ we obtain
$$
\int\limits_0^\infty h(t-\tau)(\mathscr{M}^{\gamma,2m}_\tau f)(x)d\tau=
$$
$$
=e^t\frac{2^{m-\frac{3}{2}}}{|\xi|^{m-\frac{1}{2}}}\,\Gamma\left(\frac{n+|\gamma|}{2} \right)\Gamma\left(\frac{2m-n-|\gamma|+1}{2} \right)\times
$$
$$
\times{\bf j}_\gamma(x,\xi) \int\limits_0^\infty e^{-\tau}\tau^{m-\frac{1}{2}}J_{m-\frac{1}{2}}(|\xi|\tau)d\tau
$$
$$
=e^t\frac{2^{m-\frac{3}{2}}}{|\xi|^{m-\frac{1}{2}}}\,\Gamma\left(\frac{n+|\gamma|}{2} \right)\Gamma\left(\frac{2m-n-|\gamma|+1}{2} \right)
\times
$$
$$
\times {\bf j}_\gamma(x,\xi)
\frac{2^{m-\frac{1}{2}} |\xi|^{m-\frac{1}{2}} \left(|\xi|^2+1\right)^{-m} \Gamma (m)}{\sqrt{\pi }}=
$$
$$
=\frac{\Gamma (m)\Gamma\left(\frac{n+|\gamma|}{2} \right)\Gamma\left(\frac{2m-n-|\gamma|+1}{2} \right)}{2^{2-2m}\sqrt{\pi}\left(1+|\xi|^2\right)^{m} }\,e^t  {\bf j}_\gamma(x,\xi).
$$
Let calculate the constant
$$
\frac{|S_1^+(n)|_\gamma}{N(2m,\gamma,n)}=
$$
$$
=\frac{\prod\limits_{i=1}^n{\Gamma\left(\frac{\gamma_i{+}1}{2}\right)}}{2^{n-1}\Gamma\left(\frac{n{+}|\gamma|}{2}\right)}\frac{\sqrt{\pi}}{2^{2m-n-1}\prod\limits^{n}_{i=1}\Gamma\left(\frac{\gamma_i+1}{2}\right)
	\Gamma\left(\frac{2m-n-|\gamma|+1}{2}\right)\Gamma\left(m\right)}=
$$
$$
=\frac{2^{2-2m}\sqrt{\pi}}{\Gamma\left(m\right)\Gamma\left(\frac{n{+}|\gamma|}{2}\right)\Gamma\left(\frac{2m-n-|\gamma|+1}{2}\right)}.
$$
Therefore
$$
\frac{|S_1^+(n)|_\gamma}{N(2m,\gamma,n)}\left(\frac{\partial^2}{\partial t^2}-\Delta_\gamma \right)^m \int\limits_0^\infty h(t-\tau)(\mathscr{M}^{\gamma,2m}_\tau f)(x)d\tau=
$$
\begin{equation}\label{LF}
=\frac{1}{\left(1+|\xi|^2\right)^{m}}\left(\frac{\partial^2}{\partial t^2}-\Delta_\gamma \right)^m {e^t  {\bf j}_\gamma(x,\xi)}= e^t  {\bf j}_\gamma(x,\xi),
\end{equation}
that gives $f(x)={\bf j}_\gamma(x,\xi)$. 
In the  formula (\ref{LF}) we used the fact that $\Delta_\gamma \mathbf{j}_\gamma(x;\xi)=-|\xi|\mathbf{j}_\gamma(x;\xi)$ \cite{LPSh1}.

This result confirmed by the formula (see \cite{SitnikShishkina Elsevier})
$$(M^\gamma_\rho)_x{\bf j}_\gamma(x,\xi)={\bf j}_\gamma(x,\xi)\,\,j_{\frac{n+|\gamma|}{2}-1}(\rho|\xi|).$$


\begin{thebibliography}{99.}%
	\bibitem{Finch}  Finch, D.,  Patch, S. K.,  Rakesh: Determining a function from its mean values over a family of spheres. SIAM J. Math. Anal. \textbf{35(5)}, 1213-–1240 (2004)	
	
	\bibitem{Kuchment}  Kuchment, P.: The Radon Transform and Medical Imaging.
	Society for Industrial and Applied Mathematics Philadelphia, PA, USA  (2014)
	
	\bibitem{Rubin1}
	Rubin, B.:  Fractional Integrals and Potentials. Addison-Wesley, Essex (1996)
	
	\bibitem{Rubin2} Rubin, B.:  Introduction to Radon Transforms: With Elements of Fractional Calculus and Harmonic Analysis. Cambridge University Press, UK (2015)
	
	\bibitem{Rubin3} Rubin, B.:
	Inversion formulae for the spherical mean in odd dimensions and the Euler--Poisson--Darboux equation.  Inverse Problems \textbf{24(2)},
	1--10 (2008)
	
	
	\bibitem{Agranovsky1}
	Agranovsky, M.,  Finch, D.,   Kuchment, P.: Range conditions for a spherical mean
	transform. Inverse Probl. Imaging \textbf{3(3)}, 373--382 (2009)
	
	\bibitem{Agranovsky2}  Agranovsky, M.,  Kuchment, P.,  Kunyansky, L.: On reconstruction formulas and algorithms
	for the thermoacoustic tomography. In Photoacoustic imaging and spectroscopy,
	ed. by L. Wang, CRC Press, 89--102 (2009)
	
	\bibitem{Kunyansky}  Kunyansky, L.: Explicit inversion formulae for the spherical mean Radon transform.
	Inverse Problems \textbf{23} 373--383 (2007)
	
	
	\bibitem{Weinstein2} Weinstein, A.: Spherical means in spaces of constant curvature. Annali di Matematica Pura
	ed Applicata \textbf{4(60)} 87--91 (1962)
	
	\bibitem{Dunkl}
	Hamma, M.E., Daher, R.: Estimate of K-functionals and modulus of smoothness
	constructed by generalized spherical mean operator. Pro Indian Acad. Sci. (Math. Sci.) \textbf{124(2)} 
	235--242 (2014)
	
	\bibitem{Elouadih} Elouadih, S., Daher, R.: Generalization of Titchmarsh's Theorem for the Dunkl Transform
	in the Space $L^p(\mathbb{R}^d,\omega_l(x)dx)$. International Journal of Mathematical Modelling \& Computations \textbf{6(4)} 261--267 (2016)
	
	
	\bibitem{LPSh1}  Lyakhov, L.N.,  Polovinkin, I.P.,  Shishkina, E.L.: 
	On a Kipriyanov problem for a singular ultrahyperbolic equation.
	Differ. Equ. \textbf{50(4)} 513--525 (2014)
	
	\bibitem{LPSh2} Lyakhov, L.N.,  Polovinkin, I.P.,  Shishkina, E.L.:      Formulas for the solution of the Cauchy problem for a singular wave equation with Bessel time operator. Doklady Mathematics of the Russian Academy of Sciences  \textbf{90(3)}  737--742 (2014)
	
	\bibitem{Masood}    Masood, K.,  Messaoudi, S.A., Zaman, F.D.:
	Initial inverse problem in heat equation with Bessel operator. International Journal of Heat and Mass Transfer \textbf{45(14)} 2959--2965 (2002)
	
	\bibitem{Yurko}	 Yurko, V.: 
	Inverse problems for Bessel-type differential equations on noncompact graphs using spectral data. Inverse Problems \textbf{27(4)} 045002 1--18 (2011)
	
	
	\bibitem{Krav1} Kravchenko, V.V.:  On a method for solving the inverse scattering problem on the line. Mathematical Methods in the Applied Sciences   \textbf{42(4)} 1321--1327 (2019)
	
	\bibitem{Krav2} Delgado, B.B.,  Khmelnytskaya, K.V., Kravchenko, V.V.: The transmutation operator method for efficient solution of the inverse Sturm--Liouville problem on a half-line. Mathematical Methods in the Applied Sciences.
	Mathematical Methods in the Applied Sciences  \textbf{42(18)} 7359--7366 (2019)
	
	\bibitem{KKK QuantumFest} Karapetyants, A.N.,  Khmelnytskaya, K.V., 
	Kravchenko, V.V.: A practical method for solving the inverse quantum scattering
	problem on a half line. Journal of Physics: Conference Series. \textbf{1540(1)} 012007 1--8 (2019)
	
	\bibitem{Kr2019JIIP} Kravchenko, V.V.: On a method for solving the inverse
	Sturm--Liouville problem. Journal of Inverse and Ill-posed Problems
	\textbf{27}  401--407 (2019)
	
	
	\bibitem{Zhit}    Zhitomirskii, Ya.I.:   Cauchy's problem for systems of linear partial differential equations with differential operators of Bessel type.   Mat. Sb. (N.S.) \textbf{36(78):2} 299--310 (1955)
	
	\bibitem{kipr}  Kipriyanov, I.A.:  Singular Elliptic Boundary Value Problems.
	Nauka, Moscow (1997)
	
	\bibitem{Lev1}  Levitan, B.M.: Expansion in Fourier series and integrals with Bessel functions.   Uspekhi Mat. Nauk \textbf{6:2(42)}  102--143 (1951) 
	
	\bibitem{SitnikShishkina Elsevier}  Shishkina, E.L.,  Sitnik, S.M.:
	Transmutations, singular and fractional differential equations with
	applications to mathematical physics. Elsevier, Amsterdam (2020)
	
	\bibitem{ShishkinaKir} Shishkina, E.L.: 
	Inversion of the mixed Riesz hyperbolic B-potentials.
	International Journal of Applied
	Mathematics \textbf{30(6)} 487--500 (2017)
	
	\bibitem{ShishkinaRUDN}  Shishkina, E.L.:. General Euler--Poisson--Darboux equation and hyperbolic B-potentials.
	Partial differential equations, CMFD, PFUR, M. \textbf{65(2)}  157--338 (2019)
	
	
	
	\bibitem{Nogin3}  Nogin, V.A.,   Sukhinin, E.V.: Inversion and characterization of hyperbolic potentials in $L_p$-spaces.
	Dokl. Acad. Nauk \textbf{329(5)}  550--552 (1993)
	
	
	
	\bibitem{kurant}
	Courant R., Hilbert D.: Equations of mathematical physics. T.1: Scientific ed. (1951)
	
	\bibitem{IR2}  Prudnikov, A.P.,   Brychkov, Yu.A.,  Marichev, O.I.:   Integrals and Series, Vol. 2,  Special Functions.
	Gordon \& Breach Sci. Publ., New York (1990)
	
	
	
	
	\bibitem{ShishkinaSitnik} Shishkina, E.L., Sitnik, S.M.:  General form of the Euler--Poisson--Darboux equation and application of the transmutation method. Electron. J. Differential Equations \textbf{177} 1--20 (2017)
	
	
	
\end{thebibliography}
\end{document}